\documentclass[12pt]{article}
\usepackage{amsmath,amssymb,amsthm}

\newtheorem{theorem}{Theorem}

\newtheorem{corollary}[theorem]{Corollary}

\def\barr{\begin{array}}
\def\earr{\end{array}}

\title{Minimal non-Iwasawa finite groups}
\author{Marius T\u arn\u auceanu}
\date{February 11, 2018}

\begin{document}

\maketitle

\begin{abstract}
    In this note, we describe first the structure of minimal non-Iwasawa finite groups. Then we determine the
    minimal non-Iwasawa finite groups which are modular. Also, we find connections between minimal non-Iwasawa 
    finite groups and the subgroup commutativity degree, and we give an example of a family of non-nilpotent
    modular finite groups $G_n$, $n\in\mathbb{N}$, whose subgroup commutativity degree tends to $1$ as $n$ tends to infinity.
\end{abstract}

{\small
\noindent
{\bf MSC 2010\,:} Primary 20D15; Secondary 20D30, 20D60.

\noindent
{\bf Key words\,:} minimal non-Iwasawa group, minimal non-modular group, Schmidt group, subgroup commutativity degree.}

\section{Introduction}

Let $G$ be a finite group and $L(G)$ be the subgroup lattice of $G$. A subgroup $H$ of $G$ is called
\textit{permutable} if $HK=KH$, for all $K\in L(G)$, and \textit{modular} if it is a modular element
of $L(G)$. Clearly, any normal subgroup is permutable and, by Theorem 2.1.3 of \cite{9}, a permutable
subgroup is always modular. If all subgroups of $G$ are permutable then we say that $G$ is an \textit{Iwasawa group},
while if all subgroups of $G$ are modular (that is, the lattice $L(G)$ is modular) then we say that $G$ is a
\textit{modular group}. The connection between these two classes of groups is very powerful: a finite group
$G$ is an Iwasawa group if and only if it is a nilpotent modular group (see e.g. Exercise 3, page 87, \cite{9}).
Note that a complete description of the structure of finite Iwasawa groups can be obtained by using Theorems 2.4.13
and 2.4.14 of Schmidt's book.

All groups considered in this paper are finite.

Given a class of groups $\mathcal{X}$, a group $G$ is said to be a \textit{minimal non-$\mathcal{X}$ group}, or an
\textit{$\mathcal{X}$-critical group}, if $G\notin\mathcal{X}$, but all proper subgroups of $G$ belong to $\mathcal{X}$.
Many results have been obtained on minimal non-$\mathcal{X}$ groups, for various choices of $\mathcal{X}$. For instance,
minimal non-abelian groups were analysed by Miller and Moreno \cite{6}, while Schmidt \cite{8} (see also \cite{3,7}) studied
minimal non-nilpotent groups. The latter are now known as \textit{Schmidt groups}, and their structure is as follows:
a Schmidt group $G$ is a solvable group of order $p^mq^n$ (where $p$ and $q$ are different primes) with a unique Sylow
$p$-subgroup $P$ and a cyclic Sylow $q$-subgroup $Q$, and hence $G$ is a semidirect product of $P$ by $Q$. Moreover, we have:
\begin{itemize}
\item[-] if $Q=\langle y\rangle$ then $y^q\in Z(G)$;
\item[-] $Z(G)=\Phi(G)=\Phi(P)\times\langle y^q\rangle$, $G'=P$, $P'=(G')'=\Phi(P)$;
\item[-] $|P/P'|=p^r$, where $r$ is the order of $p$ modulo $q$;
\item[-] if $P$ is abelian, then $P$ is an elementary abelian $p$-group of order $p^r$ and $P$ is a minimal normal subgroup of $G$;
\item[-] if $P$ is non-abelian, then $Z(P)=P'=\Phi(P)$ and $|P/Z(P)|=p^r$.
\end{itemize}We also recall the class of minimal non-modular $p$-groups, whose structure has been investigated in \cite{4}.
\bigskip

Our main result is the following.

\begin{theorem}
A finite group is a minimal non-Iwasawa group if and only if it is either a minimal non-modular $p$-group
or a Schmidt group $G=PQ$ with $P$ modular.
\end{theorem}

We observe that $A_4$ and $SL(2,3)$ are examples of Schmidt groups as in Theorem 1 with $P\cong\mathbb{Z}_2^2$ abelian
and $P\cong Q_8$ non-abelian, respectively.
\bigskip

By using Theorem 1, we are able to determine those minimal non-Iwasawa groups which are modular.

\begin{corollary}
The class $\cal C$ of minimal non-Iwasawa groups which are modular consists of all Schmidt groups $G=PQ$ with $P$ cyclic of order $p$.
\end{corollary}

Note that an interesting subclass of $\cal C$ is constituted by the non-trivial semidirect products $\mathbb{Z}_3\rtimes\mathbb{Z}_{2^n}$, $n\in\mathbb{N}^*$. Also, from the proof of Corollary 2, it will follow that:

\begin{corollary}
$\cal C$ is contained in the class of minimal non-cyclic groups.
\end{corollary}

\section{Proofs of the main results}

\noindent\textbf{Proof of Theorem 1.} Obviously, a minimal non-modular $p$-group is also a minimal non-Iwasawa group. Let $G=PQ$ be a Schmidt group with $P$ modular. Then $G$ is non-nilpotent, and consequently non-Iwasawa. Since all proper subgroups of $G$ are nilpotent, it suffices to prove that they are modular. Also, we may restrict to maximal subgroups. By looking to the structure of $G$ described in Section 1, we infer that these are $P\times\langle y^q\rangle$ and $\Phi(P)Q_i$, $i=1,2,...,n_q$, where $Q_i$, $i=1,2,...,n_q$, denote the conjugates of $Q$. Being a direct product of modular groups of coprime orders, $P\times\langle y^q\rangle$ is modular. On the other hand, for each $i$ we have $$\Phi(P)Q_i/Z(\Phi(P)Q_i)=\Phi(P)Q_i/\Phi(P)\langle y^q\rangle\cong Q_i/\langle y^q\rangle\cong\mathbb{Z}_q,$$implying that $\Phi(P)Q_i$ is abelian. Therefore all maximal subgroups of $G$ are modular, as desired.

Conversely, assume that $G$ is a minimal non-Iwasawa group. We distinguish the following two cases.
\smallskip

\hspace{5mm} {\bf Case 1.} $G$ is nilpotent.\\
\noindent Then $G$ is not modular. Let $G=\prod_{i=1}^k G_i$ be the decomposition of $G$ as a direct product of Sylow subgroups. Then we must have $k=1$ because all proper subgroups of $G$ are modular. Thus $G$ is a minimal non-modular $p$-group.
\smallskip

\hspace{5mm} {\bf Case 2.} $G$ is not nilpotent.\\
\noindent Then $G$ is a Schmidt group, say $G=PQ$ with $P$ and $Q$ as we described in Section 1. Since $P$ is a proper subgroup of $G$, it must be modular by our assumption. This completes the proof.
\hfill\rule{1,5mm}{1,5mm}
\bigskip

\noindent\textbf{Proof of Corollary 2.} Let $G$ be a group contained in $\cal C$. By Theorem 1, it follows that $G$ is a Schmidt of type $G=PQ$ with $P$ modular. Since $G$ is modular, so is $G_1=G/Z(G)$. But $G_1$ is again a Schmidt group of order $p^rq$ which can be written as semidirect product of an elementary abelian $p$-group $P_1$ of order $p^r$ by a cyclic group $Q_1$ of order $q$. Suppose that $r>1$. Then $P_1$ contains a proper non-trivial subgroup, say $P_2$. It is easy to see that the subgroups $$Z(G), P_1, P_2, Q_1, \mbox{ and } G_1$$form a pentagon in $L(G_1)$, a contradiction. So, $r=1$. This leads to $|P/P'|=p$, i.e. $P$ is abelian, and consequently $|P|=p$.
\hfill\rule{1,5mm}{1,5mm}
\bigskip

\noindent\textbf{Proof of Corollary 3.} Let $G$ be a group contained in $\cal C$. From the proof of Theorem 1 it follows that all maximal subgroups of $G$ are cyclic. Hence $G$ is a minimal non-cyclic group.
\hfill\rule{1,5mm}{1,5mm}

\section{Minimal non-Iwasawa finite groups and subgroup commutativity degrees}

A notion strongly connected with Iwasawa groups is the subgroup commutativity degree of a finite group $G$, defined in \cite{10} by
$$sd(G)=\frac{1}{|L(G)|^2}\,|\lbrace (H,K)\in L(G)^2 \ | \ HK=KH\rbrace|.$$This measures the probability that two subgroups of $G$ commute, or e\-qui\-va\-len\-tly that the product of two subgroups is again a subgroup. Clearly, we have $sd(G)=1$ if and only if $G$ is an Iwasawa group. $sd(G)$ has been generalized to the relative subgroup commutativity degree of a subgroup $H$ of $G$ (see \cite{11}): $$sd(H,G)=\frac{1}{|L(H)||L(G)|}\,|\lbrace (H_1, G_1) \in L(H)\times L(G) \ | \ H_1G_1=G_1H_1\rbrace|.$$These notions lead to two functions on $L(G)$, namely
$$f,g:L(G)\longrightarrow [0,1], f(H)=sd(H) \mbox{ and } g(H)=sd(H,G),\, \forall\, H\in L(G),$$whose study is proposed in \cite{11}. We remark that they are constant on each conjugacy class of subgroups of $G$. On the other hand, we have $$|Im\,f|=1\Leftrightarrow|Im\,g|=1\Leftrightarrow G = \mbox{Iwasawa group}.$$Having in mind these results, it is an interesting problem to determine the classes ${\cal C}_f$ and ${\cal C}_g$ of finite groups $G$ such that $|Im\,f|=2$ and $|Im\,g|=2$, respectively. We mention that ${\cal C}_g$ has been studied in \cite{5}. Also, it is clear that the minimal non-Iwasawa groups are contained in ${\cal C}_f$.

Another interesting problem concerning the subgroup commutativity degree is to find some natural families of groups $G_n$, $n\in\mathbb{N}$, whose subgroup commutativity degree tends to a constant $a\in[0,1]$ as $n$ tends to infinity. For $a=0$ many examples of such families are known (see e.g. \cite{1,2,10,12}). Since minimal non-Iwasawa groups have many commuting subgroups, we expect that they will have large subgroup commutativity degrees. Indeed, this is confirmed for groups in the class $\cal C$, as shows our following theorem.

\begin{theorem}
Let $p$ and $q$ be two primes such that $p\equiv 1\, (mod\, q)$, and $G_n$ be a group of order $pq^n$ contained in $\cal C$. Then $$\lim_{n\to\infty}sd(G_n)=1.$$
\end{theorem}

\noindent\textbf{Proof.} Under the notation in the previous sections, we easily infer that $L(G_n)$ consists of $G_n$, of all conjugates $Q_i$, $i=1,2,...,n_q=p$, of $Q$, and of all subgroups contained in $P\times\langle y^q\rangle\cong\mathbb{Z}_{pq^{n-1}}$. Then $$|L(G_n)|=2n+1+p.$$Also, we observe that the non-normal subgroup of $G_n$ are $Q_i$, $i=1,2,...,p$. For each $H\leq G_n$, let $C(H)=\lbrace K\in L(G) \ | \ HK=KH\rbrace$. Then $$C(H)=L(G), \,\forall\, H\unlhd\, G$$and $$C(Q_i)=L(G)\setminus\{Q_1,...,Q_{i-1},Q_{i+1},...,Q_p\}, \,\forall\, i=1,2,...,p.$$One obtains
$$sd(G_n)=\frac{1}{|L(G)|^2}\sum_{H\leq G}\!|C(H)|=\frac{1}{|L(G)|^2}\left(\sum_{H\unlhd\, G}\!|C(H)|+\sum_{H\ntrianglelefteq\, G}\!|C(H)|\right)\vspace{-3mm}$$
\begin{align*}
&\hspace{-8mm}=\frac{1}{|L(G)|^2}\left[(|L(G)|-p)|L(G)|+p(|L(G)|-p+1)\right]\\
&\hspace{-8mm}=1-\frac{p^2-p}{(2n+1+p)^2}\,,
\end{align*}which clearly tends to $1$ as $n$ tends to infinity.
\hfill\rule{1,5mm}{1,5mm}

\section{Further research}

We end our note by indicating two natural open problems concerning the above results.
\bigskip\newpage

\noindent{\bf Problem 1.} Determine the finite groups $G$ containing a unique non-Iwasawa proper subgroup $H$.

\noindent\hspace{5mm} Note that in this case $H$ must be a minimal non-Iwasawa group, and also a characteristic maximal subgroup of $G$. 
Several examples of such groups are the quasi-dihedral group $QD_{16}=\langle x,y \ | \ x^8=y^2=1,\ yxy=x^3\rangle$ and the direct products $S_3\times\mathbb{Z}_p$, where $p$ is an odd prime.
\bigskip

\noindent{\bf Problem 2.} Give a complete description of the class ${\cal C}_f$.

\noindent\hspace{5mm} Note that ${\cal C}_f$ contains any direct product between a minimal non-Iwasawa group and an Iwasawa group of coprime orders. Moreover, if $G$ is a group in ${\cal C}_f$, then by choosing a subgroup $H$ of $G$ which is minimal with the property $1\neq sd(H)=sd(G)$, it follows that $H$ is minimal non-Iwasawa. So, ${\cal C}_f$ is strongly connected to the class of minimal non-Iwasawa groups.

\vspace*{3ex}\small

\hfill
\begin{minipage}[t]{5cm}
Marius T\u arn\u auceanu \\
Faculty of  Mathematics \\
``Al.I. Cuza'' University \\
Ia\c si, Romania \\
e-mail: {\tt tarnauc@uaic.ro}
\end{minipage}

\end{document}